\documentclass[final,leqno,onefignum,onetabnum]{siamltex1213}
\usepackage{amsmath}
\usepackage{amssymb}
\usepackage{mathrsfs}
\usepackage{amsfonts}
\usepackage{multirow}
\usepackage{float}
\usepackage{enumerate}
\usepackage{bm}
\usepackage{bbm}
\usepackage{algorithm}
\usepackage{algorithmicx}
\usepackage{algpseudocode}
\usepackage[justification=centering]{caption}
\newcommand{\di}{\mathrm{d}}
\newcommand{\pal}{\partial}
\newcommand{\I}{\mathcal{I}}
\newcommand{\F}{\mathcal{F}}
\newcommand{\C}{\mathcal{C}}
\newcommand{\G}{\mathcal{G}}
\newcommand{\Y}{\mathcal{Y}}
\newcommand{\X}{\mathcal{X}}
\newcommand{\R}{\mathbb{R}}
\newcommand{\Ey}{\mathcal{E}_y}
\newcommand{\Eyw}{\mathcal{E}_{yw}}
\renewcommand{\O}{\mathcal{O}}
\renewcommand{\L}{\mathcal{L}}
\renewcommand{\P}{\mathbb{P}}
\newcommand{\INT}{\mathbb{I}}
\newcommand{\QUA}{\mathbb{Q}}
\renewcommand{\aa}{\bm{a}}
\newcommand{\bb}{\bm{b}}
\newcommand{\ii}{\bm{i}}
\newcommand{\jj}{\bm{j}}
\newcommand{\kk}{\bm{k}}

\newcommand{\xx}{\bm{x}}
\newcommand{\yy}{\bm{y}}
\newcommand{\zz}{\bm{z}}
\newcommand{\ff}{\bm{f}}
\newcommand{\uu}{\bm{u}}
\newcommand{\xxi}{\bm{\xi}}
\newcommand{\ssi}{\bm{\sigma}}
\newcommand{\bbe}{\bm{\beta}}

\newcommand{\abs}[1]{\left\vert#1\right\vert}
\newcommand{\E}[1]{\mathbb{E}\left[#1\right]}
\newcommand{\Ec}[3]{\mathbb{E}_{#1}^{#2}\left[#3\right]}

\newcommand{\CE}[2]{\mathbb{E}\left[\left.#1\right\vert #2\right]}

\newtheorem{remark}[theorem]{Remark}

\title{Efficient spectral sparse grid approximations for solving multi-dimensional forward backward SDEs\thanks{This work was supported by the National Natural Science Foundations of China under grants
91530118 and 11571351}}

\author{Yu Fu\footnotemark[2]
\and Weidong Zhao\footnotemark[3]
\and Tao Zhou\footnotemark[4]}

\begin{document}
\maketitle

\renewcommand{\thefootnote}{\fnsymbol{footnote}}

\footnotetext[2]{School of Mathematics \& Finance Institute, Shandong University, Jinan 250100, China (nielf0614@126.com)}
\footnotetext[3]{School of Mathematics \& Finance Institute, Shandong University, Jinan 250100, China (wdzhao@sdu.edu.cn)}
\footnotetext[4]{LSEC, Institute of Computational Mathematics, Academy of Mathematics and Systems Science,
Chinese Academy of Sciences, Beijing 100190, China (tzhou@lsec.cc.ac.cn)}

\renewcommand{\thefootnote}{\arabic{footnote}}

\slugger{mms}{xxxx}{xx}{x}{x--x}

\begin{abstract}
This is the second part in a series of papers on multi-step schemes for solving coupled forward backward stochastic differential equations (FBSDEs). We extend the basic idea in our former paper [W. Zhao, Y. Fu and T. Zhou, SIAM J. Sci. Comput., 36 (2014), pp. A1731-A1751] to solve high-dimensional FBSDEs, by using the spectral sparse grid approximations. The main issue for solving high dimensional FBSDEs is to build an efficient spatial discretization, and deal with the related high dimensional conditional expectations and interpolations. In this work, we propose the sparse grid spatial discretization. We use the sparse grid Gaussian-Hermite quadrature rule to approximate the conditional expectations. And for the associated high dimensional interpolations, we adopt an spectral expansion of functions in polynomial spaces with respect to the spatial variables, and use the sparse grid approximations to recover the expansion coefficients. The FFT algorithm is used to speed up the recovery procedure, and the entire algorithm admits efficient and high accurate approximations in high-dimensions, provided that the solutions are sufficiently smooth. Several numerical examples are presented to demonstrate the efficiency of the proposed methods.
\end{abstract}

\begin{keywords}
Sparse grid approximations, forward backward stochastic differential equations, conditional expectations.
\end{keywords}

\begin{AMS}
60H35, 65C20, 60H10
\end{AMS}

\pagestyle{myheadings}
\thispagestyle{plain}
\markboth{Fu, Zhao, and Zhou}{Sparse grid methods for FBSDEs}

\section{Introduction}

Backward stochastic differential equation (BSDE) in the linear sense was first introduced
by J.M.Bismut in 1973 \cite{Bismut}. Then, in 1990, Pardoux and Peng showed the existence
and uniqueness of the adapted solution for nonlinear BSDE for the first time \cite{PP1990}.
Since then, forward backward stochastic differential equations(FBSDEs) have been extensively studied,
and have been shown disperse applications in different fields, such as stochastic optimal control, nonlinear filtering,
nonlinear expectations, ect.
The general high-dimensional FBSDEs defined on a complete probability space $(\Omega,\F,\mathbb{F},\P)$ take the following form
\begin{equation}\label{FBSDEs}
\left\{\begin{aligned}
X_t &= X_0 + \int_0^t \bb(s,X_s,Y_s,Z_s)\di s+\int_0^t \ssi(s,X_s,Y_s,Z_s)\di W_s, &\text{(SDE)}\\
Y_t &= \xi + \int_t^T \ff(s,X_s,Y_s,Z_s)\di s-\int_t^T Z_s \di W_s, &\text{(BSDE)}
\end{aligned}\right.
\end{equation}
where $t\in[0,T]$ with $T$ being the fixed time horizon;
$\mathbb{F} = (\F_t)_{0\leq t\leq T}$ is the natural filtration of
the standard $d$-dimensional Brownian motion $W = (W_t)_{0\leq t\leq T}$;
$X_0\in\F_0$ and $\xi\in\F_T$ are the initial and terminal conditions for
the forward stochastic differential equation(SDE) and BSDE respectively;
$\bb:[0,T]\times\R^q\times\R^p\times\R^{p\times d}\rightarrow\R^q$ is called the drift coefficient, while
$\ssi:[0,T]\times\R^q\times\R^p\times\R^{p\times d}\rightarrow\R^{q\times d}$
is referred to the diffusion coefficient;
$\ff:[0.T]\times\R^q\times\R^p\times\R^{p\times d}\rightarrow\R^p$ is the generator of BSDE;
$(X_t,Y_t,Z_t):[0,T]\times \Omega \rightarrow \R^q\times\R^p\times\R^{p\times d}$ are the unknowns.
It is worth to note that $\bb(\cdot,\xx,\yy,\zz)$, $\ssi(\cdot,\xx,\yy,\zz)$ and $\ff(\cdot,\xx,\yy,\zz)$ are all
$\F_t$-adapted for fixed $\xx$, $\yy$ and $\zz$,
and that the two stochastic integrals with respect to $W_s$ are of the It\^o type.
A triple $(X_t,Y_t,Z_t)$ is called an $L^2$-adapted solution for
FBSDEs \eqref{FBSDEs} if it is $\F_t$-adapted, square integrable and satisfies the FBSDEs \eqref{FBSDEs}.
FBSDEs \eqref{FBSDEs} is called \emph{decoupled} when $\bb$ and $\ssi$ are both independent
of $Y$ and $Z$. In this paper, we consider the numerical solution for FBSDEs \eqref{FBSDEs}
with $\xxi=\varphi(X_T)$, $\bb(t,X_t,Y_t,Z_t)$, $\ssi(t,X_t,Y_t,Z_t)$ and $\ff(t,X_t,Y_t,Z_t)$
being deterministic functions.

Due to the complex solution structure, solutions of FBSDEs in closed form can seldom be constructed.
However, for decoupled FBSDEs, Peng \cite{Peng} introduced the following nonlinear Feynman-Kac formula, which established a deep relationship between parabolic PDEs and FBSDEs \eqref{FBSDEs}: consider the following parabolic PDE
\begin{equation}\label{PDE}
\left\{\begin{aligned}
&\uu_t(t,\xx) + \L \uu(t,\xx) + \ff(t, \xx, \uu(t,\xx), \uu_{\xx}\ssi(t,\xx)) = 0, \quad t\in[0,T), x\in\R^q\\
&\uu(T, \xx) = \varphi(\xx),
\end{aligned}\right.
\end{equation}
where
\begin{equation}
\L \uu(t,\xx) = \sum_{i=1}^{q}b_i\pal_{x_i}\uu(t,\xx) + \sum_{i,j=1}^{q}\left[\ssi\ssi^\top\right]_{i,j}\pal_{x_ix_j}\uu(t,\xx).
\end{equation}

If the above PDE \eqref{PDE} has a classical solution $\uu(t,\xx)\in\mathcal{C}^{1,2}$, then the triple $\big(X_t,\uu(t,X_t),\uu_x(t,X_t)\ssi(t,X_t)\big)$
solves the decoupled FBSDEs \eqref{FBSDEs}.

From the numerical point of view, one can use the above connections between PDEs and FBSDEs to design the so called probabilistic numerical methods for PDEs, by solving the equivalent FBSDEs. While there are a lot of works dealing with numerical schemes for BSDEs \cite{BT,BZ,JF,Cubature,LGW,Z1,ZZJ,ZCP,RO}, however, there are only a few work on numerical methods for FBSDEs \cite{DMP,MT,MPJ,MSZ,TZZ2015,ZFZ2014,ZZwJ} and the second order FBSDEs (which are related to fully non-linear PDEs) \cite{FTW,GZZ,KZZ,KZZ1}.

The key issue in numerical methods is to balance the accuracy and computational complexity. Typically, the computational complexity increases dramatically as the dimension increases. Some of the above mentioned works are designed with high order accuracy that can however only be used to deal with low dimensional FBSDEs. While some of them are low order numerical methods that are suitable for solving high dimensional problems. In particular, we highlight the work \cite{LGW}, where the constructed numerical methods can deal with very high dimensional BSDEs, however, the convergence rate is only 1/2.  We also mention the work \cite{GZZ}, where a numerical example for a 12-dimensional coupled FBSDE is reported, and it is shown by numerical test that the numerical method converges with order 1.

In this work, we aim to design high order numerical schemes for multi-dimensional FBSDEs, by extending our previous work in [W. Zhao, Y. Fu and T. Zhou, SIAM J. Sci. Comput., 36 (2014), pp. A1731-A1751]. The main difficulty for solving high dimensional FBSDEs is to efficiently evaluate the high dimensional conditional expectations and design the corresponding high dimensional interpolations. Here we shall use the sparse grid approximation technique to deal with these issues, and the FFT algorithm will be used to speed up the interpolation procedure. Several multi-dimensional examples with dimension up to 6 will be presented, and high order convergence rates up to 3 will be shown.

The rest of the paper is organized as follows. In Section 2, we present some
preliminaries, and review the multi-step schemes introduced in our previous work \cite{ZFZ2014}. The sparse grid approximation will be discussed in Section 3, and this is followed by our fully discrete numerical schemes for high dimensional FBSDEs in Section 4. In Section 5, we shall present several numerical examples.
Finally in Section 6 we give some concluding remarks. Now we introduce some notations to be used.

\begin{enumerate}[\indent\it 1.]
\item $A^\top$: the transpose of a vector/matrix $A$.
\item $\abs{\cdot}$: the Euclidean norm in the Euclidean space $\R$, $\R^q$ and $\R^{q\times d}$.
\item $\mathcal F_s^{t,\xx}$: the $\sigma$-algebra generated by the diffusion process
      $\{X_r,t\leq r\leq s, X_t=\xx\}$.
\item $\Ec{s}{t,\xx}{\eta}$: the conditional expectation of the random variable $\eta$
      under $\mathcal F_s^{t,\xx}$, i.e.,
      $\Ec{s}{t,\xx}{\eta} = \CE{\eta}{\mathcal F_s^{t,\xx}},$ and we denote $\Ec{t}{t,\xx}{\eta}$ by $\Ec{t}{\xx}{\eta}$
      for simplicity.
\item The symbols in bold denote the corresponding vector or multi-index.
      For a multi-index $\ii =(i_1,i_2,\ldots,i_d)$, $|\ii|_1 \triangleq \sum_{k=1}^d i_k$.
\item $\mathbb{F}_q^k(D)$ denotes the set of functions defined by
      \[
      \left\{\ff:D\rightarrow\R\Big\vert\mathcal{D}^{\bm\alpha} \ff\text{ continuous if }
      \alpha_i\leq k \text{ for all } i\right\},
      \]
      where $\bm{\alpha}=(\alpha_1,\alpha_2,\ldots,\alpha_q)$, $D\subset\R^q$ is
      a bounded domain and
      \[
      \mathcal{D}^{\bm{\alpha}} \ff = \frac{\pal^{|\bm{\alpha}|_1} \ff}
      {\pal x_1^{\alpha_1}\cdots\pal x_q^{\alpha_q}}
      \]
\item Let $\{\phi_i(x)\}_{i=1}^N$ be a sequence of functions from $\R$ to $\R$.
      For a multi-index $\ii=(i_1,i_2,\ldots,i_q)$,
      $\phi_{\ii}(\xx):\R^q\to \R$ is defined by
      \[
      \phi_{\ii}(\xx) = \prod_{k=1}^q \phi_{i_k}(x_k),
      \]
      where $\xx=(x_1,x_2,\ldots,x_q)$.
\end{enumerate}

\section{Time discretization for FBSDEs}

In this section, we shall briefly review the time-discrete schemes proposed in \cite{ZFZ2014}.
To this end, let us consider the uniform time partition over $[0,T],$
\[
\mathbb{T}:0=t_0<t_1<t_2<\cdots<t_N=T,
\]
with $\Delta t=\frac{T}{N}$ and $t_i=i\Delta t, i = 0,1,\ldots,N$. We set $\Delta W_{n,k}=W_{t_{n+k}}-W_{t_n}$,
$\Delta W_{t_n,t} = W_t-W_{t_n}$, $\Delta t_{n,k}=t_{n+k}-t_k$ and $\Delta t_{t_n,t}=t-t_n$ for convenience.
Then under certain regularity assumptions, for $\xx\in \R^q$ the following two reference equations
can be derived from the FBSDEs \eqref{FBSDEs} (for details, one can refer to \cite{ZFZ2014}):
\begin{eqnarray}
\frac{\di \Ec{t_n}{\xx}{Y_t}}{\di t} &=& -\Ec{t_n}{\xx}{f(t,X_t,Y_t,Z_t)},\label{dif_y} \\
\frac{\di \Ec{t_n}{\xx}{Y_t(\Delta W_{t_n,t})^\top}}{\di t} &=& -\Ec{t_n}{\xx}{f(t,X_t,Y_t,Z_t)
(\Delta W_{t_n,t})^\top}
+ \Ec{t_n}{\xx}{Z_t}. \label{dif_z}
\end{eqnarray}
The multi-step schemes in \cite{ZFZ2014} rely on efficiently approximating the derivatives in \eqref{dif_y} and \eqref{dif_z}. The following classical approximations are used for functions $u(t)\in C_b^{k+1}:$
\begin{equation}
u(t_{n+i}) = \sum_{j=0}^k\frac{(i\Delta t)^j}{j!}\frac{\di^j u}{\di t^j}(t_n) + \O\left(\Delta t_{n,i}\right)^{k+1},
\quad i = 0, 1,\ldots, k.
\end{equation}
One can easily check that $\frac{\di u}{\di t}(t_n)$ admits the following representation
\begin{equation}
\frac{\di u}{\di t}(t_n) = \sum_{i=0}^k\alpha_{k,i} u(t_{n+i})
+ \O(\sum_{i=0}^k\alpha_{k,i} (i\Delta t)^{k+1}),
\end{equation}
where $\alpha_{k,i}, i=0,1,\ldots,k$ solve the following linear system
\begin{equation}\label{alpha_equ}
\left[
\begin{array}{ccccc}
1&1&1&\cdots&1\\
0&1&2&\cdots&k\\
0&1^2&2^2&\cdots&k^2\\
\vdots&\vdots&\vdots&\vdots&\vdots\\
0&1^k&2^k&\cdots&k^k
\end{array}
\right]\cdot
\left[
\begin{array}{c}
\alpha_{k,0}\Delta t\\
\alpha_{k,1}\Delta t\\
\alpha_{k,2}\Delta t\\
\vdots\\
\alpha_{k,k}\Delta t
\end{array}
\right]=
\left[
\begin{array}{c}
0\\
1\\
0\\
\vdots\\
0
\end{array}
\right].
\end{equation}
Now by inserting the similar approach to equations \eqref{dif_y} and \eqref{dif_z} one gets
\begin{flalign}
\sum_{i=0}^k\alpha_{k,i}\Ec{t_n}{\xx}{Y_{t_{n+i}}} = -f(t_n,\xx,Y_{t_n},Z_{t_n}) + R_{\yy,n}^{k},\label{ref_y}\\
\sum_{i=1}^k\alpha_{k,i}\Ec{t_n}{\xx}{Y_{t_{n+i}}(\Delta W_{n,i})^\top} = Z_{t_n} + R_{\zz,n}^k.\label{ref_z}
\end{flalign}
To propose the multi-step schemes in \cite{ZFZ2014}, one needs the following local property of the generator of diffusion processes (see also in \cite{ZFZ2014} for details)
\begin{theorem}
Let $t_0<t$ be a fixed time, and $\xx_0\in \R^q$ be a fixed space point. If $\ff\in C^{1,2}([0,T]\times\R^q)$
and $\Ec{t_0}{\xx_0}{\abs{\L \ff(t,X_t)}}<+\infty$, then we have
\begin{equation}\label{gnrt_eq1}
\left. \frac{\di\Ec{t_0}{\xx_0}{f(t,X_t)}}{\di t} \right|_{t=t_0} = \left.\frac{\di\Ec{t_0}{\xx_0}{\ff(t,\bar X_t)}}{\di t}\right|_{t=t_0},
\end{equation}
where $\bar X_t$ is a diffusion process satisfying
\begin{equation}
\bar{X}_t = \xx + \int_{t_0}^t\bar{\bb}_s\di s + \int_{t_0}^t\bar\ssi_s\di W_s
\end{equation}
with
$\bar \bb_s=\bar \bb(s,\bar X_s; t_0,\xx_0),$ $\bar \ssi_s = \bar{\ssi}(s,\bar X_s; t_0,\xx_0)$
being smooth functions of $(s,\bar X_s)$ with parameters $(t_0,\xx_0)$ that satisfy
\begin{equation*}
\bar \bb(t_0,\bar X_{t_0}; t_0,\xx_0)= \bb(t_0,\xx_0),\quad \bar{\ssi}(t_0,\bar X_{t_0}; t_0,\xx_0)=\ssi(t_0,\xx_0).
\end{equation*}
\end{theorem}
By combining this local property and equations (2.6)-(2.7), the following multi-step numerical schemes are proposed in \cite{ZFZ2014}
\begin{algorithm}[H]
\caption{Multi-step semi-discrete schemes for coupled FBSDEs}
Assume that $Y^{N-i}$ and $Z^{N-i}$ ($i = 0,1,\ldots,k-1$) are known.
For $n=N-k,\ldots,0$, solve
$X^{n,j}(j=1,2,\dots,k)$, $Y^n=Y^n(X^n)$ and $Z^n=Z^n(X^n)$ by
\vspace{0.25cm}
\begin{enumerate}
\item  Set $Y^{n,0}(X^n)=Y^{n+1}(X^n)$ and $Z^{n,0}(X^n)=Z^{n+1}(X^n)$.
\vspace{0.25cm}
\item Set $l=0$ and let $\epsilon_0$ be a given tolerance. Solve $Y^{n,l+1}(X^n)$ and $Z^{n,l+1}(X^n)$ by the
following steps,
\begin{align*}
&X^{n,j} = X^n + \bb(t_n,X^n,Y^{n,l},Z^{n,l})\Delta t_{n,j} + \ssi(t_n,X^n,Y^{n,l},Z^{n,l})\Delta W_{n,j}, \,\, j = 1,\ldots,k, \\
&Z^{n,l+1}(X^n) = \sum_{j=1}^k\alpha_{k,j}\Ec{t_n}{X^n}{\bar Y^{n+j}(\Delta W_{n,j})^\top},\\
&\alpha_{k,0}Y^{n,l+1}(X^n) = -\sum_{j=1}^k\alpha_{k,j}\Ec{t_n}{X^n}{\bar Y^{n+j}}
-\ff(t_n,X^n,Y^{n,l+1},Z^{n,l+1}),
\end{align*}
until
\[
\max\Big\{|Y^{n,l+1}(X^n)-Y^{n,l}(X^n)|,|Z^{n,l+1}(X^n)-Z^{n,l}(X^n)|\Big\}< \epsilon_0,
\]
where $\bar Y^{n+j}$ are the values of $Y^{n+j}$ at the space point $X^{n,j}$.
\vspace{0.25cm}
\item Let $Y^n(X^n) = Y^{n,l+1}(X^n)$ and $Z^n(X^n)=Z^{n,l+1}(X^n)$.
\end{enumerate}
\end{algorithm}

The above schemes are semi-discrete schemes proposed in \cite{ZFZ2014}, and it was shown that the $k$-th step scheme admits a $k$ order convergence rate, provided that $1\leq k \leq 6.$  To efficiently solve the FBSDEs, however, we also have to introduce a spatial discretizition, and should guarantee a high quality spatial approximation, e.g., the approximation of conditional expectations, to balance the entire numerical error. In fact, this is the main purpose of this work, and we will build efficient algorithms to deal with this issue in the following sections. Meanwhile, it is noted that an iterative procedure is used in the above schemes, this is due to the couple property of FBSDEs. For Decoupled FBSDEs, such a procedure may be omitted.

\section{Sparse discretization and corresponding function approximations}
In the last section, we have introduced the high order semi-discrete schemes for coupled FBSDEs. However, to make the multi-step schemes more efficient, one should design efficient numerical methods for evaluating the conditional expectations and high dimensional interpolations. In \cite{ZFZ2014}, one uses the uniform tensor spatial meshes, and the tensorized Gaussian Hermite quadrature rule was used to evaluate the conditional expectations, moreover, the Lagrange interpolation method is used to compute the non-grid information. It was shown that such a combination is less efficient for high dimensional FBSDEs, as the required computational work increases exponentially as the dimension increases. To this end, we shall introduce the sparse grid approximation method, which is introduced originally for approximating high dimensional integral \cite{BNR,Smolyak}, to deal with this issue. We remark that the sparse grid approximation has been used in many different research topics, see e.g., \cite{NTW,XH,SY1,SY2,ZGZ} and references therein.

\subsection{The multi-dimensional sparse grids}

In this section, we follow closely the idea and notations in \cite{Smolyak} to given a basic introduction for constructing high dimensional sparse grids.  To begin, let us consider a sequence of grids $\{\chi^i\}_{i=1}^\infty$
in $\R$
\[
\chi^i = \{x_0^i, x_1^i, \ldots, x_{N_i-1}^i\}, \quad i = 1, 2, \ldots.
\]
The total number $N_i$ of points is usually chosen to be $N_i=2^i+1.$ Based on such an one dimensional sequence, one can build the $q$-dimensional sparse grids via
\begin{equation}\label{spg}
\chi_q^p = \bigcup_{q\leq |\ii|_1\leq p} \chi^{i_1}\otimes\chi^{i_2}\otimes\cdots\otimes\chi^{i_q},
\end{equation}
where $p\geq q$ is an integer, and $\ii=(i_1,...,i_q)$ is a multi-index. We call the one-dimensional sequence \emph{nested} if it satisfies $\chi^1\subset\chi^2\cdots\subset\chi^k\subset\cdots.$
In such cases, one can rearrange the sequence into the following hierarchical order:
\[
\chi = \chi^1\cup(\chi^2\setminus\chi^1)\cup\cdots\cup(\chi^k\setminus\chi^{k-1})
\cup\cdots = \{x_0,x_1,\ldots,x_{N_{k}-1},\ldots\}.
\]
Let $\I^i=\{0,1,\ldots,N_i-1\}$ with $\I^0={\emptyset}$, then $\chi^i=\{x_j, j\in\I^i\}$.
The \emph{nested} structure admits many advantages in constructing high dimensional sparse grids. From the computational cost point of view, one needs only count each different point once, and then \eqref{spg} can be written as
\begin{equation}\label{nested_spg}
\chi_q^p = \bigcup_{q\leq |\ii|_1\leq p}\tilde\chi^{i_1}\otimes\tilde\chi^{i_2}\otimes
\cdots\otimes\tilde\chi^{i_q}\quad \text{with } \quad \tilde\chi^i\triangleq
\left\{\begin{array}{ll}
\chi^1, & i=1,\\
\chi^i\setminus\chi^{i-1}, & i>1 .
\end{array}\right.
\end{equation}
In this paper, we shall adopt two papular types of one dimensional grids to construct the high dimensional sparse grid approximations. The first one is the sparse grid based on
the one dimensional Chebyshev-Gauss-Lobatto(CGL) grids $\{\C^i\}_{i=1}^{\infty}$ that is defined by
$$\C^i = \big\{x_j^i = \cos(\frac{j\pi}{2^i}), \,\,\, j=0,\ldots, 2^i\big\}.$$
It is easy to see that the CGL sequence is nested. 
Using such a nested property, in our framework, the sparse grid based on CGL points will be used to build high dimensional function approximations and related high dimensional interpolations.

The other type of sparse grid we shall use is the one based on the one dimensional Gauss-Hermite(GH) grid $\{\G^i\}_{i=1}^{\infty}$ which is defined as
$$\G^i = \{x_j^i, j=1,\ldots,2^i-1\},$$
where $\{x_j^i, j=1,\ldots,2^i-1\}$ are roots of the Hermite polynomial
of order $2^i-1.$  Obviously, the GH sequence $\{\G^i\}_{i=1}^{\infty}$ is not nested. Note that the above definitions are for the one dimensional case, and one can construct the associated high dimensional sparse grids ($\G_q^p$ and $\C_q^p$) based on these one-dimensional sets and the formula (\ref{spg}).  One can refer to Fig.1 to have a first glance at the two types of sparse grids in two dimensions. In this work, we shall use the sparse grids quadrature rule based on the GH sequence to approximate the associated high dimensional integrals (conditional expectations) in our multi-step schemes for solving FBSDEs.
\begin{figure}[H]
\centering
\includegraphics[width=2.5in]{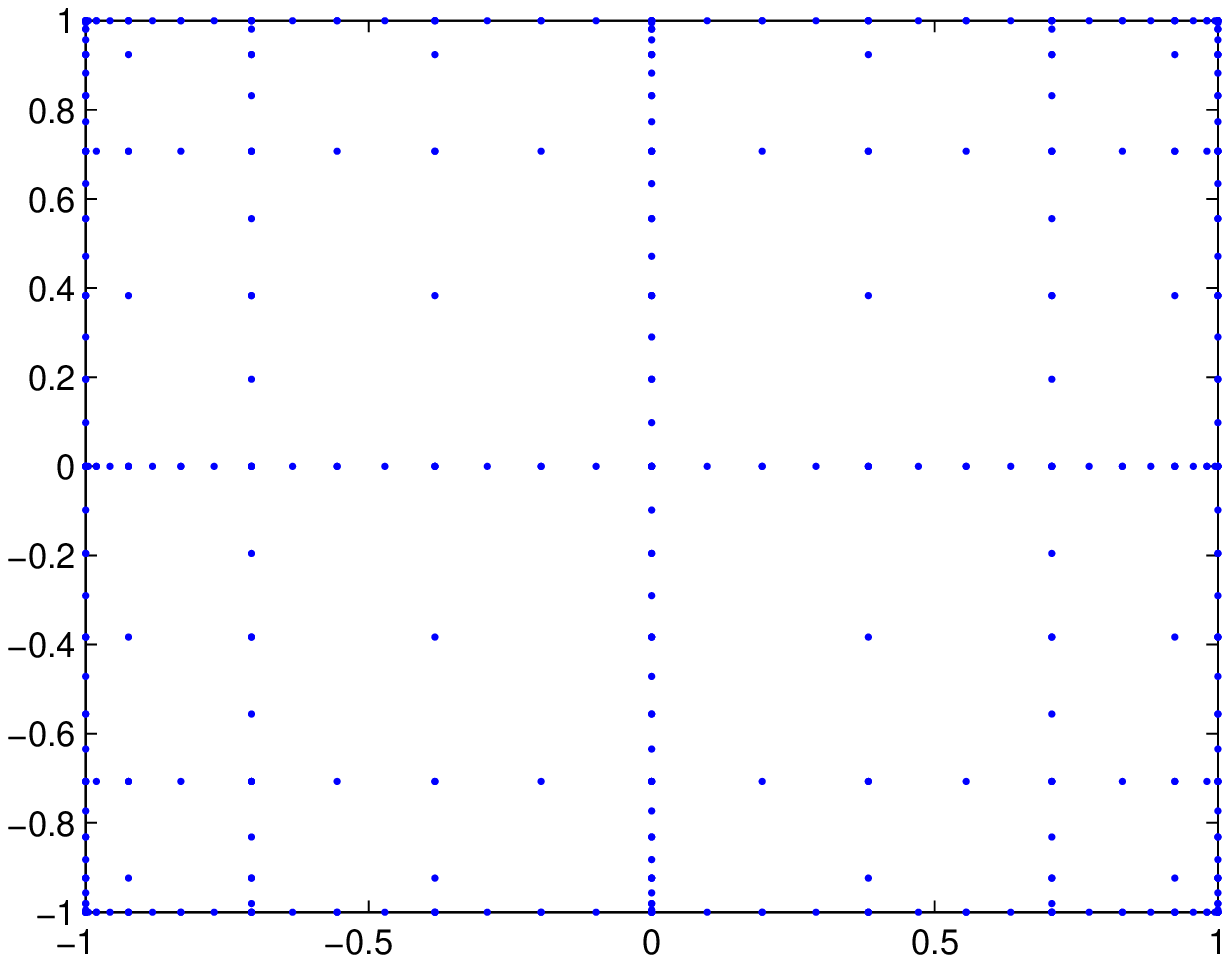}
\includegraphics[width=2.5in]{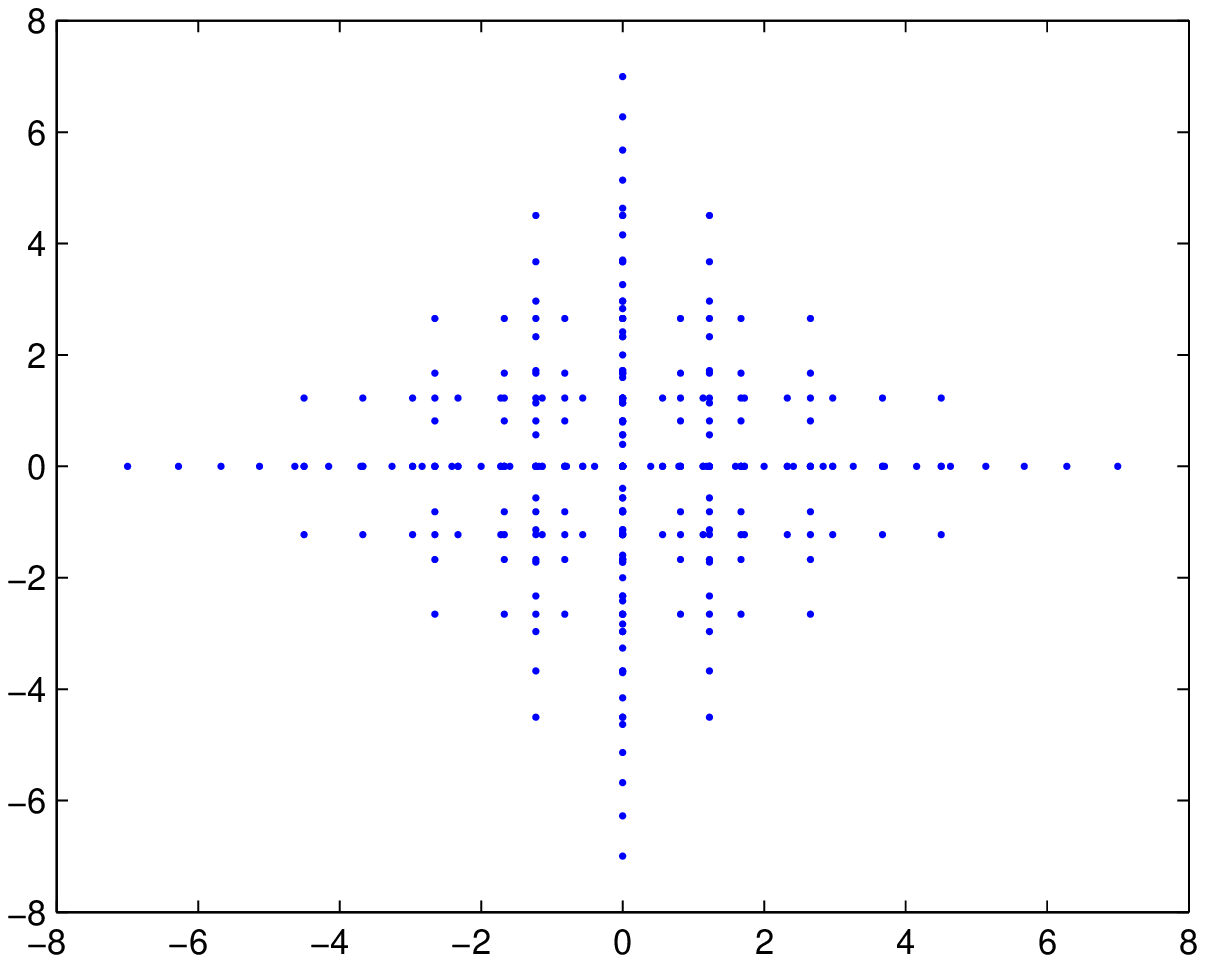}
\caption{Sparse grid for CGL $\C_2^6$ (Left) and sparse grid for GH $\G_2^6$ (Right).}
\end{figure}

\subsection{Function approximations on sparse grids}
In this section, we discuss function approximations on sparse grids, and this will play an important role in our numerical schemes for solving FBSDEs. To this end, let us begin with the one-dimensional case,
and consider a bounded interval $I\subset\R$ in which we have a set of $N_i$ points
$\chi^i = \{x_k, k\in\I^i\}.$ Let $\omega(x)>0$ be a weight function in $I$
and $\{\phi_k\}_{k\in\I^i}$ be a set of basis functions (usually orthogonal polynomials) in $L_{\omega}^2(I).$
Then for $f\in C(I)$, if the values of $f$ at $x\in\chi^i$ are known, we can
construct an interpolation approximation of $f$ via
\begin{eqnarray*}
\mathbb{I}^if(x) &=& \sum_{k\in\I^i}b_k^i\phi_k(x),
\end{eqnarray*}
where the coefficients $\{b_k^i\}$ are obtained by solving the following linear system
\begin{eqnarray*}
f(x_j^i)=\sum_{k\in\I^i}b_k^i\phi_k(x_j^i), \quad j\in\I^i.
\end{eqnarray*}
Then, one can define the multi-variate interpolation approximation for $\R^q$-functions by the following sparse interpolation operator:
\begin{equation}\label{spg_int}
\mathbb{I}_q^p [\ff] = \sum_{q\leq |\ii|_1 \leq p} \Delta^{i_1}\otimes\Delta^{i_2}
\otimes\cdots\otimes\Delta^{i_q} [\ff] \quad \text{with} \quad
\Delta^{i}\triangleq
\left\{\begin{array}{ll}
\INT^1, & i=1,\\
\INT^i-\INT^{i-1}, & i>1 .
\end{array}\right.
\end{equation}
or equivalently,
\begin{equation}\label{spg_int2}
\mathbf{\INT}_q^p[\ff] = \sum_{p-q<|\ii|_1\leq p}(-1)^{p-|\ii|_1}
\left(\begin{array}{c}
q-1\\
p-|\ii|_1
\end{array}\right)\INT^{i_1}\otimes\INT^{i_2}\otimes\cdots\otimes\INT^{i_q} [\ff],
\end{equation}
where $p\geq q$ is a positive integer, and $\INT^{i_1}\otimes\INT^{i_2}\otimes\cdots\otimes\INT^{i_q}$
stands for a $q$-dimensional interpolation based on the sparse grid information
$(\xx, f(\xx))$, $\xx\in\chi^{i_1}\otimes\chi^{i_2}\otimes\cdots\otimes\chi^{i_q}$
with the basis functions $\phi_{\kk}=\prod_{j=1}^{q}\phi_{k_j}$, $k_j\in\I^{i_j},$ and $\kk=(k_1,...,k_q).$

To further reduce the computational complexity, we shall use hierarchical bases in this work. We recall the following definition:
\begin{definition}\label{hie_bases}
For nested grids, a set of basis functions
$\{\tilde\phi(x)_k\}_{k=0}^{\infty}$ is called hierarchical,
if
\[
\tilde\phi_k(x_j)=0 \,\text{ for any }\, j\in\I^i,~k\notin \I^i.
\]
\end{definition}
By the above definition, it is easy to see that the expansion
coefficients $\{b^i_k\}_{k\in\I^i}$ under the hierarchical bases will independent of the level
index $i.$ That is, for any $f(x)\in C(I)$, we can write
\[
\INT^i f(x) = \sum_{k\in\I^i}b_k\tilde\phi_k(x),\quad
\Delta^i f(x) = \sum_{k\in\tilde\I^i}b_k\tilde\phi_k(x), \quad i=1,2,\ldots,
\]
where $\tilde\I^i=\I^i\setminus\I^{i-1}$, and the coefficients $\{b_k\}_{k\in\I^i}$ can be determined by
\[
f(x_j) = \sum_{k\in\I^i}b_k\tilde\phi_k(x_j), \quad
j\in\I^i, \quad  i = 1, 2, \ldots.
\]
Thus, in case a nested sparse grid and the corresponding hierarchical bases
are used, the interpolation procedure \eqref{spg_int} can be simplified as
\begin{equation}
\begin{aligned}
\mathbb{I}_q^p [f](\xx) =~& \sum_{q\leq |\ii|_1 \leq p} \Delta^{i_1}\otimes\Delta^{i_2}
\otimes\cdots\otimes\Delta^{i_q}[f](\xx)\\
=~& \sum_{q\leq |\ii|_1 \leq p} \sum_{\kk\in\tilde\I^{i_1}\times\cdots\times\tilde\I^{i_q}}
b_{k_1}\cdots b_{k_q}\tilde\phi_{k_1}(x_1)\cdots\tilde\phi_{k_q}(x_q)\\
=~& \sum_{\kk\in\I_q^p}b_{\kk}\tilde\phi_{\kk}(\xx),
\end{aligned}
\end{equation}
where $\xx=(x_1,\ldots,x_q)$, $\tilde\phi_{\kk}=\prod_{i=1}^{q}\tilde\phi_{k_i},$ and
\[
\I_q^p = \bigcup_{q\leq|\ii|_1\leq p} \tilde\I^{i_1}\times\cdots\times\tilde\I^{i_q}.
\]
Then, the expansion coefficients $\{b_{\kk},\kk\in\I_q^p\}$ can be solved by
\begin{equation}\label{cond_tran}
f(\xx_{\jj}) = \sum_{\kk\in\I_q^p} b_{\kk}\tilde\phi_{\kk}(\xx_{\jj}),
\quad \forall \jj\in\I_q^p.
\end{equation}
Thus,
$\INT_d^q$ defines a unique interpolation operator onto the finite space
\[
V_q^p \triangleq \textmd{span}\{\tilde\phi_{\kk}, \kk\in\I_q^p\}.
\]
In \cite{SY2}, the authors developed a fast transform to solve the equation \eqref{cond_tran}.
To illustrate the solving procedure, let's take $\I_q^p = \I_2^4$ as an example, i.e.,
\[
f(\xx_{\jj}) = \sum_{\kk\in\I_q^p} b_{\kk}\tilde\phi_{\kk}(\xx_{\jj}),
\quad \forall \jj\in\I_2^4.
\]
To make it more clearly, first of all we consider how to compute the values of
$b_{\kk}, \kk\in\I^1\times\I^1$ by the following equation,
\[
f(\xx_{\jj}) = \sum_{\kk\in\I_q^p} b_{\kk}\tilde\phi_{\kk}(\xx_{\jj}),
\quad \forall \jj\in\I^1\times\I^1.
\]
The coefficients can be computed in two steps.
\begin{enumerate}[1.]
\item Perform the following transform on $\{f(\xx_{\jj}),\jj\in\I^1\times\I^1\}$
along the first dimension, and get $\{b_{\jj}^\prime,\jj\in\I^1\times\I^1\}$.
\[
b_{k_1, j_2}^\prime = \sum_{j_1\in\I^1}f(x_{j_1},x_{j_2})T_{j_1,k_1}, \quad
k_1\in\I^1, j_2\in\I^1,
\]
where matrix $T$ is the inverse of $(\tilde\phi_k(x_j))_{k,j}$.
\item Perform the transform on $\{b_{\jj}^\prime,\jj\in\I^1\times\I^1\}$
along the second dimension, and get the coefficients $\{b_{\jj},\jj\in\I^1\times\I^1\}$.
\[
b_{k_1, k_2} = \sum_{j_2\in\I^1}b_{k_1,j_2}^\prime T_{j_2, k_1}, \quad
k_1\in\I^1, k_2\in\I^1,
\]
where matrix $T$ is the inverse of $(\tilde\phi_k(x_j))_{k,j}$.
\end{enumerate}
Based on this fact, we can compute the coefficients $\{b_{\kk}, \kk\in\I_2^4\}$
in the following way. Firstly, we perform one-dimensional transforms on
$\{f(\xx_{\jj}),\jj\in\I^3\times\I^1\}$, $\{f(\xx_{\jj}),\jj\in\I^2\times\tilde\I^2\}$
and $\{f(\xx_{\jj}),\jj\in\I^1\times\tilde\I^3\}$ along the first dimension, respectively. We get
$\{b_{\jj}^\prime, \jj\in\I_2^4\}$. Then apply one-dimensional transforms on
$\{b_{\jj}^\prime, \jj\in\I^1\times\I^3\}$, $\{b_{\jj}^\prime, \jj\in\tilde\I^2\times\I^2\}$
and $\{b_{\jj}^\prime, \jj\in\tilde\I^3\times\I^1\}$ along the second dimension.
We obtain all the values of $\{b_{\jj}, \jj\in\I_2^4\}$.

Now we extend this procedure to the $q$-dimensional case, and introduce the following algorithm.
For more details about the fast transform, readers may refer to the work \cite{SY2}.

\begin{algorithm}
  \caption{Fast transform on sparse grid.}
  \label{fast_tran}
  \begin{algorithmic}
    \Require $q$, $p$, $\chi_q^p$, $\{\ff(\xx_{\jj}),\jj\in\I_q^p\}$ and
             $\{\phi_{\kk}(\xx),\kk\in\I_q^p\}$
    \Ensure $\{\bb_{\kk},\kk\in\I_q^p|\ff(\xx_{\jj})=\sum_{\kk\in\I_q^p}
               \bb_{\kk}\tilde\phi_{\kk}(\xx_{\jj})\}$
    \Function{FastTran}{$q$,$p$,$\{\ff(\xx_{\jj}),\jj\in\I_q^p\}$}
    \State $\bb_{\jj} \gets \ff(\xx_{\jj})$, for all $\jj\in\I_q^p$
    \For{$q^\prime = 1 \to q$}
        \ForAll{$\ii^\prime\in\{\ii^\prime=(i_1,\ldots,i_{q^\prime-1},i_{q^\prime+1},\ldots,i_q)
                \big{|} |\ii^\prime|_1 \leq p-1$\}}
        \State $(T_{k,j})_{k,j\in\I^{p-|\ii^\prime|_1}} \gets
                (\tilde\phi_k(x_j))_{k,j\in\I^{p-|\ii^\prime|_1}}^{-1}$
        \State $\{\bb_{\jj},\jj\in\tilde\I^{i_1}\times\cdots\times\tilde\I^{i_{q^\prime-1}}
        \times\I^{p-|\ii^\prime|_1}
        \times\tilde\I^{i_{q^\prime+1}}\times$
        \State $\quad \cdots\times\tilde\I^{i_q}\}\gets\sum_{k\in\I^{p-|\ii^\prime|_1}}
               \bb_{j_1,\ldots,j_{d^\prime-1},
               k,j_{q^\prime+1},\ldots,j_q} T_{k,j_{q^{\prime}}}$
        \EndFor
    \EndFor
    \State \Return{$\{\bb_{\kk},\kk\in\I_q^p\}$}
    \EndFunction
  \end{algorithmic}
\end{algorithm}

%

\subsection{The approximation of conditional expectations}\label{approx_ce}
Recall our multi-step schemes in Algorithm 1 for solving FBSDEs, to eventually solve the FBSDEs, one needs to do spatial discretizations.  In this work, we aim at solving the solution pair $(Y(\xx),Z(\xx))$ on a given sparse grid $\xx \in \chi_q^p,$ on each time level. Note that other required information of $(Y,Z)$ can be obtained by using the interpolation procedure based on the sparse grid information. Although the FBSDEs are essentially defined on the unbounded domain, however, we are usually interested in the values of $(Y_0,Z_0)$
in a bounded domain of $\xx$, which we denote by $[\aa_0, \bb_0]$. Then for each time level $t_n$,
we only need to approximate the solution of $(Y_{t_n},Z_{t_n})$ in a bounded domain $[\aa_n, \bb_n]$
so that the values of $(Y_{t_n},Z_{t_n})$ outside will not influence our computation of $(Y_0,Z_0)$
in $[\aa_0, \bb_0]$. More precisely, for the time level $t=t_n,$ we aim at solving $(Y_{t_n}(\xx),Z_{t_n}(\xx))$ for $\xx \in \chi_q^{p_n} \subset [\aa_n, \bb_n].$ To do this, one has to approximate the associated conditional expectations in Algorithm 1. As mentioned before, we shall adopt the sparse grid GH quadrature rule to approximate the conditional expectations.
To this end, we propose the following choice for the bounded domain: for any $\xx\in \chi_q^{p_n},$ there holds
\begin{equation}\label{eq:domain}
\xx + \bb(t_n,\xx,\yy,\zz)\cdot j\Delta t \pm \ssi(t_n,\xx,\yy,\zz)\cdot
\sqrt{2j\Delta t}M \mathbbm{1}_{[d\times1]}\in[\aa_{n+j},\bb_{n+j}]
\end{equation}
for every $j=1,2,\ldots,k,$  where $M$ is defined by
$$M=\max\{ |\zz|_\infty, \,\, \zz\in  \G^p_q\},$$
where $\G^p_q$ is the sparse grid GH points that are used to approximate the conditional expectations,
and $\mathbbm{1}_{[d\times1]}=(1,1,\ldots,1)^{\top},$ and $|\zz|_\infty = \max |z_k|.$

The motivation for the above choice is that one needs to include all information that is used in the sparse grid GH quadrature rule. More precisely, consider our multi-step schemes in Algorithm 1, and suppose that (for $j=1,2,\ldots,k$) the values of $Y^{n+j}(\xx)$ on $\xx  \in \chi_q^{p_{n+j}}$ are known. That is, we consider the $k$-step scheme.
Then, we can compute the value of $Y^{n+j}(\xx)$ for every $\xx\in[\aa_{n+j}, \bb_{n+j}]$ by
\begin{equation}\label{eq:interpolation}
\INT_q^{p_{n+j}}\big[Y^{n+j}(\xx)\big] = \sum_{\ii\in\I_q^{p_{n+j}}}
\bbe_{\ii}^{n+j}\phi_{\ii}(\xx), \quad j=1,2,\ldots,k,
\end{equation}
where the coefficients $\{\bbe_{\ii}\}_{\ii\in\I_q^{p_{n+j}}}$ can be obtained by Algorithm \ref{fast_tran} based on the sparse grid information $\xx  \in \chi_q^{p_{n+j}},$ and $\INT_q^{p_{n+j}}$ is the associated interpolation operator.

Now, to solve $(Y^n, Z^n)$ for $\xx  \in \chi_q^{p_n},$ we need to approximate the associated conditional expectations, e.g., $\Ec{t_n}{\xx}{\bar Y^{n+j}}.$  Note that we have
\begin{align*}
\Ec{t_n}{\xx}{\bar Y^{n+j}} &=\E{Y^{n+j}(\xx+\bb(t_n,\xx,\yy,\zz)\cdot j\Delta t
+ \ssi(t_n,\xx,\yy,\zz)\Delta W_{n,j})}\\
& =\int_{\R^d} \Big(Y^{n+j}(\xx+\bb(t_n,\xx,\yy,\zz)\cdot j\Delta t +
\ssi(t_n,\xx,\yy,\zz)\cdot \xxi) \\
&\qquad \cdot\frac{1}{\sqrt{(2\pi j\Delta t)^q}}\exp{\left\{-\frac{\xxi^\top\xxi}{2j\Delta t}\right\}}
\Big)\di \xxi\\
&=\frac{1}{\sqrt{\pi^q}}\int_{\R^q}\Big(Y^{n+j}\big(\xx+\bb(t_n,\xx,\yy,\zz)\cdot j\Delta t +
\ssi(t_n,\xx,\yy,\zz)\cdot\sqrt{2j\Delta t}\xxi\big)\\
&\qquad\qquad\cdot\exp\big\{-\xxi^\top\xxi\big\}\Big) \di \xxi.
\end{align*}
In \cite{ZFZ2014}, we have proposed the tensor grid of GH points to approximation the above integral, namely,
\[\begin{aligned}
\Ec{t_n}{\xx}{\bar Y^{n+j}} &\approx \QUA_q^{p,T} \left[Y^{n+j}(\bar X^{n+j})\right] \\
&= \sum_{\ii=1}^L\omega_{\ii}Y^{n+j}\Big(\xx+\bb(t_n,\xx,\yy,\zz)\cdot j\Delta t +
\ssi(t_n,\xx,\yy,\zz)\cdot\sqrt{2j\Delta t}\xxi_{\ii}\Big),
\end{aligned}\]
where $\xxi_{\ii}$ and  $\omega_{\ii}$ are the GH
quadrature points and the corresponding weights, respectively. We have defined the associated quadrature operator by $\QUA_q^{p,T}.$
Then, it is clear to see that one needs the following function information
\begin{align}
\zz=\xx+b(t_n,\xx,\yy,\zz)\cdot j\Delta t +
\sigma(t_n,\xx,\yy,\zz)\cdot\sqrt{2j\Delta t}\xxi_{\ii}.
\end{align}
Thus, we can see that the choice of \eqref{eq:domain} is reasonable because it makes every point $z$
stays inside the interval $[\aa_{n+j}, \bb_{n+j}]$.

However, when the dimension of the Brownian motion grows higher, the computational
requirements of the above tensor quadrature rule increase exponentially with respect to the dimension.
Hence we shall resort to the sparse grid GH quadrature rule, in this work, to approximate the conditional expectations. We denote by $\QUA_q^p$ the sparse grid Gauss-Hermite quadrature operator, namely,
\begin{equation}\label{spg_qua}
\QUA_q^p = \sum_{p-q<|\ii|_1\leq p}(-1)^{p-|\ii|_1}
\left(\begin{array}{c}
q-1\\
p-|\ii|_1
\end{array}\right)\QUA^{i_1}\otimes\QUA^{i_2}\otimes\cdots\otimes\QUA^{i_q},
\end{equation}
where the quadrature operator $\QUA^k$ denotes the one-dimensional Gauss-Hermite
quadrature based on the grid $\G^k$. Then, the conditional expectations involved in Algorithm 1 can be approximated by
\[\begin{aligned}
&\Ec{t_n}{\xx}{\bar Y^{n+j}} = \QUA_q^{p_{n+j}}\Big[Y^{n+j}\big(\xx+\bb(t_n,\xx,\yy,\zz)\cdot j\Delta t +
\ssi(t_n,\xx,\yy,\zz)\cdot\sqrt{2j\Delta t}\xxi\big)\Big] + \mathcal{R}_{\mathbb{Q}}^{y,n+j}, \\
&\Ec{t_n}{\xx}{\bar Y^{n+j}\big(\Delta W_{n,j})^\top} = \QUA_q^{p_{n+j}}\Big[Y^{n+j}(\xx+
\bb(t_n,\xx,\yy,\zz)\cdot j\Delta t +
\ssi(t_n,\xx,\yy,\zz)\cdot\sqrt{2j\Delta t}\xxi\big)\\
&\qquad\quad\cdot \sqrt{2j\Delta t}\xxi^\top\Big] + \mathcal{R}_{\mathbb{Q}}^{yw,n+j}, \quad j = 1,2,\ldots,k,
\end{aligned}\]
where the two terms $\mathcal{R}_{\mathbb{Q}}^{y,n+j}$ and $\mathcal{R}_{\mathbb{Q}}^{yw,n+j}$
are quadrature errors.

As discussed in \cite{BNR,NR}, for functions $f\in\mathbb{F}_1^k([-1,1])$, the error of the one-dimensional
Clenshaw-Curtis quadrature rule is given by,
\[
\abs{\int_{-1}^1 f(x)\di x - \bar\QUA_1^p[f]} \leq C\mathcal{N}_Q^{-k},
\]
where $\bar\QUA$ stands for the quadrature using Clenshaw-Curtis points,
$\mathcal{N}_Q$ is the number of quadrature points and the constant
$C$ relies on the upper bound of the $k$-th derivative of $f$.
From this one-dimensional estimates,
the authors obtain the following result of the sparse grid quadrature
for functions defined on a high-dimensional cube $f\in\mathbb{F}_q^k([-1,1]^q)$,
\[
\abs{\int_Df(x)\di x - \bar\QUA_q^p[f]}
\leq C_{q,k}\cdot\mathcal{N}_Q^{-k}\cdot\big(\log(\mathcal{N}_Q)\big)^{(k+1)(q-1)},
\]
where $\mathcal{N}_Q$ is the number of spare grid quadrature points and the constant
$C_{q,k}$ depends on $q$ and the upper bound of the $k$-th derivative of $f$.
Analogously for functions $f\in\mathbb{F}_1^k(\R)$, we have the following error estimates
in \cite{SSO},
\[
\int_\R f(x)e^{-x^2}\di x - \QUA_1^p[f] \leq C\mathcal{N}_Q^{-k/2}.
\]
Therefore, by conducting a similar procedure in \cite{BNR,NR}, we can obtain the following theorem.
\begin{theorem}\label{err_qua}
For any functions $f\in\mathbb{F}_q^k(\R^q)$, the error of sparse grid
quadrature rule $\QUA_q^p$ is
\[
\abs{\int_{\R^q}\ff(\xx)\exp\{-\xx^\top\xx\}\di \xx - \QUA_q^p[\ff]}
\leq C_{q,k}\cdot\mathcal{N}_Q^{-k/2}\cdot\big(\log(\mathcal{N}_Q)\big)^{(k/2+1)(q-1)},
\]
where $\mathcal{N}_Q = \mathcal{N}_Q(q,p)$ is the number of the sparse grid quadrature
points used by $\QUA_q^p$
and the constant $C_{q,k}$ depends only on $q$ and the upper bound of
the $k$-th derivative of $\ff$.
\end{theorem}

Note that the above quadrature rule uses non-grids information (i.e., not all grid points $\xx$ belong to $\chi_q^{p_{n+j}}$. As mentioned, these information will be obtained via the sparse grid interpolation, i.e., the procedure (\ref{eq:interpolation}). To this end, we propose the fully discrete schemes for approximating the conditional expectations:
\begin{align}
\Ec{t_n}{\xx}{\bar Y^{n+j}}&= \QUA_q^{p_{n+j}}\Big[\INT_q^{\tilde{p}_{n+j}}\big[
Y^{n+j}(\xx+\bb(t_n,\xx,\yy,\zz)\cdot j\Delta t \nonumber \\
&\qquad + \ssi(t_n,\xx,\yy,\zz)\cdot\sqrt{2j\Delta t}\xxi)\big]\Big]
+ \mathcal{R}_{\mathbb{QI}}^{y,n+j} + \mathcal{R}_{\mathbb{Q}}^{y,n+j}, \label{eq:inter_quadra_1} \\
\Ec{t_n}{\xx}{\bar Y^{n+j}(\Delta W_{n,j})^\top} &= \QUA_q^{p_{n+j}}\Big[
\INT_q^{\tilde{p}_{n+j}}\big[Y^{n+j}(\xx+
\bb(t_n,\xx,\yy,\zz)\cdot j\Delta t \nonumber\\
&\qquad + \ssi(t_n,\xx,\yy,\zz)\cdot\sqrt{2j\Delta t}\xxi)\big]
\cdot\sqrt{2j\Delta t}\xxi^\top\Big]
+ \mathcal{R}_{\mathbb{QI}}^{yw,n+j} + \mathcal{R}_{\mathbb{Q}}^{yw,n+j},\label{eq:inter_quadra_2}
\end{align}
where $\tilde{p}_{n+j}, p_{n+j}\geq q$ are positive integers, and the interpolation errors are defined by
\[
\mathcal{R}_{\mathbb{QI}}^{y,n+j} = \QUA_q^{\tilde{p}_{n+j}}\Big[Y^{n+j}-\INT_q^{\tilde{p}_{n+j}}[Y^{n+j}]\Big],
\quad \mathcal{R}_{\mathbb{QI}}^{yw,j}
= \QUA_q^p\Big[\big(Y^{n+j}-\INT_q^{\tilde{p}_{n+j}}[Y^{n+j}]\big)\cdot\sqrt{2j\Delta t}\xxi^\top\Big].
\]
We note that both the interpolation error and the quadrature error can be well controlled, provided that the functions admit certain regularities, and meanwhile, suitable sparse grids are used. For detailed error estimates of sparse grid quadratures \& interpolations, one can refer to \cite{BNR,NTW,SY1}.
\begin{remark}
In the above discussions, we have used the notations $\QUA_q^{p_{n+j}}$ and $\INT_q^{\tilde{p}_{n+j}}$ to stand for the quadrature and interpolation operators. In both notations, we use $\tilde{p}_{n+j}$ and $p_{n+j}$ to specify the different levels of sparse grids that are used for each time level $t=t_{n+j}.$ We remark that one can of course use uniform sparse grid in each time level, however, one would benefit if different level can be used according to the different time level, from the view of computational cost. In particular, for the time level with a larger domain $[\aa_k,\bb_k],$ one may need a high level sparse grids to obtain a good accuracy.
\end{remark}

In what follows, We are aimed to show how to find the effective computational domains $[\aa_n, \bb_n]$ for each time level $t=t_n$, and this is also a key issue  for solving FBSDEs. In the beginning of this section, we have supposed that the following
holds
\begin{equation}
\xx + \bb(t_n,\xx,\yy,\zz)\cdot j\Delta t \pm \ssi(t_n,\xx,\yy,\zz)\cdot
\sqrt{2j\Delta t}M \mathbbm{1}_{[d\times1]}\in[\aa_{n+j},\bb_{n+j}]
\end{equation}
for any $\xx\in \chi_q^{p_n},$ where $\chi_q^{p_n}$ is the set of sparse grids for $t=t_n.$ For simple cases, the coefficients of the forward SDE are bounded
functions, i.e.,
\[
\|\bb\|_\infty \leq C_{\bb},
\quad \abs{\ssi}_\infty \leq C_{\ssi}.
\]
With $[\aa_0,\bb_0]$ given, we can simply set the bounded domains $\{[\aa_{n}, \bb_{n}]\}_{n=1}^{N}$
by
\[
\begin{aligned}
\aa_{n+1} &= \aa_{n} + C_{\bb} \cdot \Delta t
- \bm{C}_{\ssi} \cdot \sqrt{2\Delta t}M \mathbbm{1}_{[d\times1]},\\
\bb_{n+1} &= \bb_{n} + C_{\bb} \cdot \Delta t
+ \bm{C}_{\ssi} \cdot \sqrt{2\Delta t}M \mathbbm{1}_{[d\times1]}.
\end{aligned}
\]
However, if $\bb$ and $\ssi$ are unbounded functions, the
determination of these intervals becomes much more complex, as the values of $\bb$ and $\ssi$ keep changing
during the iterative procedures. Nevertheless, one can always find a large computational domain to fix this issue (yet with huge computational cost). Another solution is to adopt an efficient approximation method for the whole space $\R^q$ instead of using bounded domains approximations. And this is our ongoing project.

\section{Fully discrete multi-step schemes for multi-dimensional FBSDEs}
We summarize in this section the entire multi-step schemes for solving multi-dimensional FBSDEs.
Let us consider the $k$-step scheme, and assume that we have obtained the initial values, i.e., $Y^{N-k}(\xx)$ and $Z^{N-k}(\xx)$, $i = 0,1,\ldots,k-1$, are known
for $\xx\in\C_q^{p_{N-i}} \subset [\aa_{N-i},\bb_{N-i}]$. Then, to solve $\big(Y^{n}(\xx), Z^{n}(\xx)\big)^{N-k}_{0}$, we do the following steps:
\begin{itemize}
\item We choose a suitable computational domain $[\aa_n, \bb_n]$ for $t=t_{n},$  and construct the corresponding sparse grids $ \C_q^{p_n} \subset [\aa_n, \bb_n].$ In our setting, for each time level, we construct the sparse grids by transforming the standard CGL sparse grids from $[-1,1]^q$ to $[\aa_n, \bb_n].$

\item For each $\xx\in \C_q^{p_n}$ we use Algorithm 1 to solve $\big(Y^{n}(\xx), Z^{n}(\xx)\big)$ with $X_n=\xx.$ In this procedure, we will use the approximation methods in (\ref{eq:inter_quadra_1})-(\ref{eq:inter_quadra_2}) to deal with the high dimensional conditional expectations. In particular, we shall use the following hierarchical basis functions $\{\tilde T_{\kk}(\xx)\}_{\kk\in\I_q^{p_n}}$ that is introduced in \cite{SY1}:
 \begin{align*}
&\tilde T_{\kk}(\xx) = \prod_{i=1}^{q}\tilde T_{k_i}(x_i),
\quad \text{with} \quad  \tilde T_{k_i}(x) =
\left\{\begin{array}{ll}
T_{k_i}(x), & k\in\I^1,\\
T_{k_i}(x)-T_{2^j-k_i}(x), & k_i\in\tilde\I^j,
\end{array}\right.\\
& \qquad \qquad T_{k_i}(x) = \cos\left(k_i\arccos\left(\frac{x-\bb_n^i}{\bb_n^i-\aa_n^i}
+\frac{x-\aa_n^i}{\bb_n^i-\aa_n^i}\right)\right),
\end{align*}
 i.e,  $T_k(x)$ is the transformed classical $k$-th order one dimensional Chebyshev polynomial. We note that by using this type of bases and the CGL sparse grid, the FFT algorithm can be used to speed up  the recovery procedure in Algorithm 2, see also in \cite{SY2}.

\item One do the above procedure until the $t=0.$
\end{itemize}
A detailed description for our multi-step schemes is also shown in Algorithm \ref{alg_msch}.
\begin{algorithm}
  \caption{$k$-step scheme for solving high-dimensional FBSDEs}
  \label{alg_msch}
  \begin{algorithmic}
    \Require $Y^{N-i}(\xx)$ and $Z^{N-i}(\xx)$, defined on $\xx\in\C_q^{p_i}[\aa_i,\bb_i]$, for $i=0,1,\ldots,k-1$.
    \Ensure  $Y^0(\xx)$ and $Z^0(\xx)$ on $\xx\in\C_q^{p_0}[\aa_0,\bb_0]$.
    \For {$i=0 \to k-1$}
    \State $\{\beta_{\jj}^{N-i},\jj\in\I_q^{p_{N-i}}\}$ $\gets$
           {\sc FastTran}($p_{N-i}$, $q$, $\{Y^{N-i}(\xx_{\jj}), \jj\in\I_q^{p_{N-i}}\}$).
    \State $\{\gamma_{\jj}^{N-i},\jj\in\I_q^{p_{N-i}}\}$ $\gets$
           {\sc FastTran}($p_{N-i}$, $q$, $\{Z^{N-i}(\xx_{\jj}), \jj\in\I_q^{p_{N-i}}\}$).
    \EndFor
    \For {$n=N-k \to 0$}
      \ForAll {$\ii\in\I_q^{p_n}$}
      \State $l=0$.
      \State $Y^{n,l}(\xx_{\ii}) \gets \sum_{\kk\in\I_q^{p_{n+1}}} \beta_{\kk}^{n+1}\tilde T_{\kk}(\xx_{\ii})$.
      \State $Z^{n,l}(\xx_{\ii}) \gets \sum_{\kk\in\I_q^{p_{n+1}}} \gamma_{\kk}^{n+1}\tilde T_{\kk}(\xx_{\ii})$.
      \Repeat
        \For {$m = 1\to k$}
          \ForAll {$W_{\jj}\in\G_d^P$}
          \State $\X_{\jj}^m \gets \xx_{\ii}+b(t_n,\xx_{\ii},Y^{n,l}(\xx_{\ii}),Z^{n,l}(\xx_{\ii}))\cdot m\Delta t$
          \State $\qquad~~+\sigma(t_n,\xx_{\ii},Y^{n,l}(\xx_{\ii}),Z^{n,l}(\xx_{\ii}))\cdot \sqrt{2m\Delta t}W_{\jj}$.
          \State $\Y_{\jj}^m \gets \sum_{\kk\in\I_q^{p_{n+m}}} \beta_{\kk}^{n+m}\tilde T_{\kk}(\X_{\jj}^m))$.
          \EndFor
        \State $\Ey^m \gets \sum_{\kk\in\I_q^{p_{n+m}}} \omega_{\jj} \Y_{\jj}^m$.
        \State $\Eyw^m \gets \sum_{\kk\in\I_q^{p_{n+m}}} \omega_{\jj} \Y_{\jj}^m \cdot \sqrt{2m\Delta t}W_{\jj}^\top$.
        \EndFor
      \State $Z^{n,l+1}(x_{\ii}) \gets \sum_{j=1}^k\alpha_{k,j}\Eyw^j$.
      \State $Y^{n,l+1}(x_{\ii})\gets\sum_{j=0}^k\alpha_{k,j}\Ey^j -
             f(t_n,x_{\ii},Y^{n,l+1}(x_{\ii}),Z^{n,l+1}(x_{\ii}))$.
      \State $l=l+1$.
      \Until {$\max\{\abs{Y^{n,l+1}(\xx_{\jj})-Y^{n,l}(\xx_{\jj})},
                \abs{Z^{n,l+1}(\xx_{\jj})-Z^{n,l}(\xx_{\jj})}\}<\epsilon_0$}.
      \If {$n\neq 0$}
      \State $\{\beta_{\jj}^{n},\jj\in\I_q^{p_n}\}$ $\gets$
             {\sc FastTran}($p_n$, $q$, $\{Y^{n}(\xx_{\jj}), \jj\in\I_q^{p_n}\}$).
      \State $\{\gamma_{\jj}^{n},\jj\in\I_q^{p_n}\}$ $\gets$
             {\sc FastTran}($p_n$, $q$, $\{Z^{n}(\xx_{\jj}), \jj\in\I_q^{p_n}\}$).
      \EndIf
      \EndFor
    \EndFor
  \end{algorithmic}
\end{algorithm}
\begin{remark}
Note that in the above algorithm, one needs to use some initial values, i.e.,
$\left\{Y^{N-i}(\xx), Z^{N-i}(\xx)\right\}$ for $i=0,1,\ldots,k-1$. This can be obtained by running some standard low order numerical schemes with small time steps, or use the Runge-Kutta methods \cite{JF} or DC methods \cite{TZZ2015} to do the initialization. In our numerical examples, we shall directly set these values to be known to avoid the initialization error.
\end{remark}

\section{Numerical experiments}

In this section, we provide several constructive numerical examples to show the efficiency of our numerical
algorithm for solving high-dimensional FBSDEs.
We shall also present a numerical comparison between the spectral
sparse grid(SSG) method and our previous approach in \cite{ZFZ2014} (where standard Lagrange interpolations on tensor grids (LTG) are used). To specify the main differences between the two approaches, we list the main techniques used in both methods in Table 1.

\begin{table}[H]
\small
\caption{Main techniques used in the SSG and LTG methods for solving FBSDEs. TP: tensor product.  SG: sparse grid. GH: Gaussian-Hermite.}
\begin{center}
\begin{tabular}{c|ccc}
\hline
Method  & Meshes  & Conditional expectations &  Approximation \& interpolation  \\
\hline
SSG     & sparse grid            & SG GH quadrature              & SG interpolation         \\
LTG     & TP uniform mesh            & TP GH quadrature             & Lagrangian\\
\hline
\end{tabular}
\end{center}
\end{table}

In what follows, we will denote by CR and RT the convergence rate and the running time, respectively.
In all our numerical tests, the numerical results, which include numerical errors,
convergence rates, and running times, are obtained on a computer with 16 Intel Xeon E5620 CPUs (2.40 GHz),
and 3.0 GB free RAM, coding in FORTRAN 95.

Example 1:  We first consider a two dimensional example, and we shall also report the numerical comparison between the SSG and LTG methods. In this two dimensional example, we set the components of $\bb$ as
\begin{equation*}
b_i(t,\xx,\yy,\zz)=\cos(4(x_i+t))/4-1, \quad i=1,2,
\end{equation*}
and we set
$\ssi$ to be a diagonal matrix with diagonal components
$$\ssi_{ii}=\cos(4(x_i+t))\sin(4(x_i+t))/4,\quad i=1, 2.$$ Furthermore, we choose a one dimensional generator function $f$ as
\begin{align*}
f(t,\xx,y,\zz) &=\frac12\sum_{i=1}^2 z_i\sin^2(4(t+x_i)) - y\sum_{i=1}^2\cos^2(4(t+x_i))\\
&\quad+\sin(4(t+x_2))\cos^2(4(t+x_1))(\sin(4(t+x_1))-1)\\
&\quad+\sin(4(t+x_1))\cos^2(4(t+x_2))(\sin(4(t+x_2))-1),
\end{align*}
so that the exact solution is
\begin{align*}
&\qquad Y_t = \sin(4(X_{t,1}+t))\sin(4(X_{t,2}+t))\\
&Z_{t,i} = \left(\prod_{k=1}^2\sin(4(X_{t,k}+t))\right)\cos^2(4(t+X_{t,i)}),\quad i = 1, 2.
\end{align*}
Due to the periodic property of this example, we can solve the FBSDEs in a fixed space domain
$[-\pi, \pi]^2$. For the SSG method, we shall use the  standard Chebyshev-Gauss-Lobatto sparse grid $\C_2^7$ for
the spatial discretization, and the Gauss-Hermite sparse grid $\G_2^3$ will be used for the high dimensional  quadratures (conditional expectations).
In the LTG method, as mentioned before, we shall use the tensor grid of uniform mesh, and the degree of Lagrangian interpolation polynomials $n$ are decided by the following formula
\begin{equation}\label{dx}
(\Delta x)^{(n+1)} = (\Delta t)^{(k+1)} \quad \Rightarrow  \quad \Delta x = (\Delta t)^{\frac{k+1}{n+1}}
\end{equation}
to balance the time discretization error and the interpolation error. For more detailed explanations of the LTG method, one can refer to \cite{ZFZ2014}. We now solve example 1 with the SSG method and the LTG method, and the errors, running times and convergence rates are shown in Table 2.
\begin{table}[H]
\small
\setlength{\belowcaptionskip}{0pt}
\renewcommand{\arraystretch}{0.8}
\caption{Numerical Errors and convergence rates for Example 1 by SSG (Top) and LTG (Bottom).}
\begin{center}
\begin{tabular}{c|c|c|c|c|c|c|c}
\hline
step number &  error    &   N=8   &   N=16   &   N=32   &  N=64 &  N=128 & CR         \\
\hline
\multirow{3}{*}{1-step}
& $\mathcal{E}_Y$ & 3.991E-02 & 2.050E-02 & 1.039E-02 & 5.232E-03 & 2.625E-03 & 0.982 \\
\cline{2-8}
& $\mathcal{E}_Z$ & 5.186E-02 & 2.649E-02 & 1.339E-02 & 6.733E-03 & 3.377E-03 & 0.986 \\
\cline{2-8}
& RT & 2.519 & 4.638 & 9.521 & 19.199 & 38.679 & \\
\cline{2-8}
\hline
\multirow{3}{*}{2-step}
& $\mathcal{E}_Y$ & 5.620E-03 & 1.456E-03 & 3.670E-04 & 9.182E-05 & 2.286E-05 & 1.987 \\
\cline{2-8}
& $\mathcal{E}_Z$ & 6.978E-03 & 1.847E-03 & 4.813E-04 & 1.225E-04 & 3.090E-05 & 1.955 \\
\cline{2-8}
& RT & 6.852 & 15.596 & 33.401 & 68.541 & 143.552 & \\
\cline{2-8}
\hline
\multirow{3}{*}{3-step}
& $\mathcal{E}_Y$ & 9.748E-04 & 1.342E-04 & 1.728E-05 & 2.264E-06 & 8.196E-07 & 2.632 \\
\cline{2-8}
& $\mathcal{E}_Z$ & 3.091E-03 & 3.850E-04 & 4.757E-05 & 6.009E-06 & 8.834E-07 & 2.955 \\
\cline{2-8}
& RT & 6.899 & 19.403 & 43.835 & 95.902 & 197.608 & \\
\cline{2-8}
\hline
\end{tabular}
\end{center}
\vspace{0.45cm}
\small
\setlength{\belowcaptionskip}{0pt}
\renewcommand{\arraystretch}{0.8}
\begin{center}
\begin{tabular}{c|c|c|c|c|c|c|c}
\hline
step number    & errors &   N=8   &   N=16   &   N=32   &  N=64 &  N=128 & CR         \\
\hline
\multirow{3}{*}{1-step} & $\mathcal{E}_Y$ & 7.204E-01 & 3.989E-01 & 1.849E-01 & 7.321E-02 & 2.871E-02 & 1.174 \\
                        \cline{2-8}
                        & $\mathcal{E}_Z$ & 2.670E-01 & 1.534E-01 & 7.409E-02 & 3.152E-02 & 1.332E-02 & 1.093 \\
                        \cline{2-8}
                        & RT & 1.047 & 4.609 & 24.120 & 141.209 & 844.974 &  \\
                        \cline{2-8}
\hline
\multirow{3}{*}{2-step} & $\mathcal{E}_Y$ & 4.873E-01 & 1.105E-01 & 2.220E-02 & 4.165E-03 & 7.708E-04 & 2.333 \\
                        \cline{2-8}
                        & $\mathcal{E}_Z$ & 1.999E-01 & 4.706E-02 & 1.174E-02 & 2.494E-03 & 4.159E-04 & 2.205 \\
                        \cline{2-8}
                        & RT & 3.861 & 12.970 & 60.661 & 363.052 & 2139.540 &  \\
                        \cline{2-8}
\hline
\multirow{3}{*}{3-step} & $\mathcal{E}_Y$ & 2.656E-01 & 2.252E-02 & 2.295E-03 & 2.247E-04 & 2.024E-05 & 3.401 \\
                        \cline{2-8}
                        & $\mathcal{E}_Z$ & 1.136E-01 & 9.841E-03 & 1.306E-03 & 1.417E-04 & 1.244E-05 & 3.243 \\
                        \cline{2-8}
                        & RT & 13.010 & 43.490 & 193.784 & 968.788 & 4929.977 &  \\
                        \cline{2-8}
\hline
\end{tabular}
\end{center}
\end{table}
From the above two tables, we immediately learn that both methods admit high order convergence rates, more precisely, the $k$-step schemes admits a $k$-order convergence rate (we only listed the numerical results for $1\leq k \leq 3$). However, the computational complexity exhibits big differences between the two proposed methods. For example, for 3-step methods with $N=128$ (the 7th column), the running time is 197s for the SSG methods vases 4930s for the LTG methods, and this is obvious due to the efficient sparse grid discretization and the efficient sparse grid interpolations.

Example 2: our next example is $q$-dimensional decoupled FBSDEs. More precisely, we set the components of $\bb$ as
\begin{align*}
b_i(t,\xx,y,\zz)=\frac{1}{q} x_i \exp{(-x_i^2)}, \quad i=1,...,q,
\end{align*}
and again, we set $\ssi$ to be a diagonal matrix with diagonal components
$$\ssi_{i,i}= \frac{1}{q}\exp{(-x_i^2)}, \quad i =1,...,q.$$
The generator function $f$ is chosen to be
\begin{align*}
f(t,x,y,z) &=-\sum_{i=1}^q x_iz_i + \frac{1}{q^2}y
 - \frac{1}{q^3}\sum_{i=1}^q\Big(x_i^2+\exp(-2x_i^2)\Big)\prod_{k=1\atop k\ne i}^q(x_k+t)\\
& \quad - \frac{1}{q}\sum_{i=1}^qx_i^2\sum_{j=1}^q\prod_{k=1 \atop k\ne i, k\ne j}^q(x_k + t).
\end{align*}
It is easy to show that the exact solution takes the following form
\begin{align*}
Y_t = \frac1q \sum_{j=1}^q \Big( X_{t,j}^2\prod_{k=1 \atop k\ne j}^q(X_{t,k}+t)\Big),
\end{align*}
\begin{align*}
Z_{t,i} = \frac{1}{q^2}\exp\left(-X_{t,i}^2\right)\sum_{j=1}^q \Big(X_{t,j}^2\prod_{k=1\atop k\ne i, k\ne j}^q(X_{t,k}+t)\Big)+ \frac{2X_{t,i}}{q^2}\exp(-X_{t,i}^2)\prod_{k=1\atop k\ne i}^q(X_{t,k}+t).
\end{align*}
We solve the above FBSDEs with $q=3,4,5,6$, and the numerical results are listed in Table \ref{ex_decouple}. Again, the proposed multi-step schemes admit high order convergence rates even for the 6-dimensional problem. To show the computational complexity of the schemes with respect to the dimension, we show in Fig.2 the growth of running time against the dimension, and it is seems that the running time grows in certain polynomial level (non-exponential).
\begin{table}
\small
\setlength{\belowcaptionskip}{0pt}
\caption{Errors and convergence rates for Example 2.}\label{ex_decouple}
\begin{center}
\begin{tabular}{c|c|c|c|c|c|c|c|c}
\hline
step number    & sparse grid &           &   N=8   &   N=16   &   N=32   &  N=64 &  N=128 & CR         \\
\hline
\multirow{3}{*}{1-step}
& \multirow{3}{*}{$q=3$, $\C_3^4$ \& $\G_3^4$} & $\mathcal{E}_Y$   & 5.717E-02 & 2.999E-02 & 1.557E-02 & 8.021E-03 & 4.104E-03 & 0.950 \\
\cline{3-9}
& \multirow{3}{*}{}         & $\mathcal{E}_Z$ & 3.129E-02 & 1.575E-02 & 7.988E-03 & 4.057E-03 & 2.057E-03 & 0.981 \\
\cline{3-9}
& \multirow{3}{*}{}         & RT & 0.107 & 0.226 & 0.389 & 0.627 & 1.268 &  \\
\cline{3-9}
\hline
\multirow{3}{*}{2-step}
& \multirow{3}{*}{$q=3$, $\C_3^4$ \& $\G_3^4$} & $\mathcal{E}_Y$ & 2.766E-03 & 7.861E-04 & 2.091E-04 & 5.396E-05 & 1.371E-05 & 1.918 \\
\cline{3-9}
& \multirow{3}{*}{}         & $\mathcal{E}_Z$ & 4.833E-03 & 1.197E-03 & 3.006E-04 & 7.588E-05 & 1.916E-05 & 1.994 \\
\cline{3-9}
& \multirow{3}{*}{}         & RT & 0.131 & 0.235 & 0.509 & 1.067 & 2.187 &  \\
\cline{3-9}
\hline
\multirow{3}{*}{3-step}
& \multirow{3}{*}{$q=3$, $\C_3^4$ \& $\G_3^4$} & $\mathcal{E}_Y$ & 1.244E-04 & 1.183E-05 & 1.061E-06 & 9.522E-08 & 8.587E-09 & 3.460 \\
\cline{3-9}
& \multirow{3}{*}{}         & $\mathcal{E}_Z$ & 3.425E-04 & 4.375E-05 & 5.575E-06 & 7.067E-07 & 8.922E-08 & 2.976 \\
\cline{3-9}
& \multirow{3}{*}{}         & RT & 0.173 & 0.370 & 0.750 & 1.478 & 3.072 &  \\
\cline{3-9}
\hline

\multirow{3}{*}{1-step}
& \multirow{3}{*}{$q=4$, $\C_4^5$ \& $\G_4^5$} & $\mathcal{E}_Y$ & 1.196E-01 & 6.208E-02 & 3.187E-02 & 1.626E-02 & 8.259E-03 & 0.965 \\
\cline{3-9}
& \multirow{3}{*}{}         & $\mathcal{E}_Z$ & 3.445E-02 & 1.735E-02 & 8.772E-03 & 4.439E-03 & 2.244E-03 & 0.985 \\
\cline{3-9}
& \multirow{3}{*}{}         & RT & 1.682 & 3.153 & 6.178 & 12.934 & 26.182 &  \\
\cline{3-9}
\hline
\multirow{3}{*}{2-step}
& \multirow{3}{*}{$q=4$, $\C_4^5$ \& $\G_4^5$} & $\mathcal{E}_Y$ & 9.114E-03 & 2.817E-03 & 7.851E-04 & 2.087E-04 & 5.426E-05 & 1.854 \\
\cline{3-9}
& \multirow{3}{*}{}         & $\mathcal{E}_Z$ & 6.769E-03 & 1.628E-03 & 4.038E-04 & 1.012E-04 & 2.547E-05 & 2.012 \\
\cline{3-9}
& \multirow{3}{*}{}         & RT & 2.513 & 4.442 & 9.606 & 20.160 & 41.504 &  \\
\cline{3-9}
\hline
\multirow{3}{*}{3-step}
& \multirow{3}{*}{$q=4$, $\C_4^5$ \& $\G_4^5$} & $\mathcal{E}_Y$ & 4.879E-04 & 7.611E-05 & 1.052E-05 & 1.378E-06 & 1.763E-07 & 2.865 \\
\cline{3-9}
& \multirow{3}{*}{}         & $\mathcal{E}_Z$ & 1.176E-03 & 1.420E-04 & 1.767E-05 & 2.215E-06 & 2.786E-07 & 3.009 \\
\cline{3-9}
& \multirow{3}{*}{}         & RT & 1.993 & 5.722 & 13.296 & 28.614 & 59.671 &  \\
\cline{3-9}
\hline

\multirow{3}{*}{1-step}
& \multirow{3}{*}{$q=5$, $\C_5^6$ \& $\G_5^6$} & $\mathcal{E}_Y$ & 2.269E-01 & 1.177E-01 & 6.023E-02 & 3.061E-02 & 1.549E-02 & 0.969 \\
\cline{3-9}
& \multirow{3}{*}{}         & $\mathcal{E}_Z$ & 4.345E-02 & 2.196E-02 & 1.110E-02 & 5.608E-03 & 2.829E-03 & 0.985 \\
\cline{3-9}
& \multirow{3}{*}{}         & RT & 23.196 & 47.566 & 99.318 & 203.587 & 411.991 &  \\
\cline{3-9}
\hline
\multirow{3}{*}{2-step}
& \multirow{3}{*}{$q=5$, $\C_5^6$ \& $\G_5^6$} & $\mathcal{E}_Y$ & 2.941E-02 & 8.713E-03 & 2.376E-03 & 6.223E-04 & 1.600E-04 & 1.885 \\
\cline{3-9}
& \multirow{3}{*}{}         & $\mathcal{E}_Z$ & 9.182E-03 & 2.216E-03 & 5.534E-04 & 1.393E-04 & 3.510E-05 & 2.005 \\
\cline{3-9}
& \multirow{3}{*}{}         & RT & 33.610 & 80.876 & 177.315 & 372.765 & 766.094 &  \\
\cline{3-9}
\hline
\multirow{3}{*}{3-step}
& \multirow{3}{*}{$q=5$, $\C_5^6$ \& $\G_5^6$} & $\mathcal{E}_Y$ & 2.549E-03 & 4.157E-04 & 6.185E-05 & 8.531E-06 & 1.127E-06 & 2.789 \\
\cline{3-9}
& \multirow{3}{*}{}         & $\mathcal{E}_Z$ & 2.364E-03 & 2.687E-04 & 3.254E-05 & 4.031E-06 & 5.037E-07 & 3.045 \\
\cline{3-9}
& \multirow{3}{*}{}         & RT & 36.049 & 105.507 & 245.986 & 530.853 & 1106.870 &  \\
\cline{3-9}
\hline

\multirow{3}{*}{1-step}
& \multirow{3}{*}{$q=6$, $\C_6^7$ \& $\G_6^7$} & $\mathcal{E}_Y$   & 4.068E-01 & 2.117E-01 & 1.084E-01 & 5.500E-02 & 2.779E-02 & 0.969 \\
\cline{3-9}
& \multirow{3}{*}{}         & $\mathcal{E}_Z$ & 5.899E-02 & 2.996E-02 & 1.516E-02 & 7.652E-03 & 3.855E-03 & 0.984 \\
\cline{3-9}
& \multirow{3}{*}{}         & RT & 368.274 & 792.758 & 1640.774 & 3346.529 & 6758.066 &  \\
\cline{3-9}
\hline
\multirow{3}{*}{2-step}
& \multirow{3}{*}{$q=6$, $\C_6^7$ \& $\G_6^7$} & $\mathcal{E}_Y$ & 7.236E-02 & 2.130E-02 & 5.783E-03 & 1.510E-03 & 3.869E-04 & 1.891 \\
\cline{3-9}
& \multirow{3}{*}{}         & $\mathcal{E}_Z$ & 1.307E-02 & 3.277E-03 & 8.382E-04 & 2.134E-04 & 5.405E-05 & 1.978 \\
\cline{3-9}
& \multirow{3}{*}{}         & RT & 1112.524 & 2813.641 & 6114.390 & 12814.473 & 26454.702 &  \\
\cline{3-9}
\hline
\multirow{3}{*}{3-step}
& \multirow{3}{*}{$q=6$, $\C_6^7$ \& $\G_6^7$} & $\mathcal{E}_Y$ & 1.037E-02 & 1.742E-03 & 2.528E-04 & 3.407E-05 & 4.433E-06 & 2.806 \\
\cline{3-9}
& \multirow{3}{*}{}         & $\mathcal{E}_Z$ & 3.924E-03 & 4.391E-04 & 5.339E-05 & 6.647E-06 & 8.328E-07 & 3.045 \\
\cline{3-9}
& \multirow{3}{*}{}         & RT & 594.060 & 1767.938 & 4110.876 & 8853.131 & 18481.311 &  \\
\cline{3-9}
\hline
\end{tabular}
\end{center}
\end{table}

Example 3:  we next consider the $q$-dimensional coupled FBSDEs. We set
\begin{equation*}
b_i(t,\xx,y,\zz) =\frac{t}{2}\cos^2(y+x_i), \quad i=1,...,q.
\end{equation*}
And the diagonal matrix $\ssi$ is chosen with components $\ssi_{i,i}=\frac{t}{2}\sin^2(y+x_i).$ The generator function is chosen as
\begin{align*}
f(t,x,y,z) &=\sum_{i=1}^q z_i - \frac{1}{q}\left(1+\frac{t}{2}\right)\sum_{i=1}^q x_i^2
- \frac{t}{q}\left(\sum_{i=1}^{q-1} x_i(x_{i+1}+t)+x_q(x_1+t)\right)\\
& \quad - \frac{t^2}{2M^2}\left(\sum_{i=1}^{q-1}(x_{i+1}+t)\sin^4(y+x_i)+(x_1+t)\sin^4(y+x_q)\right).
\end{align*}
It can be checked that the exact solution is
\begin{align*}
&\qquad \qquad \qquad \quad  Y_t = \frac1q\sum_{j=1}^q X_{t,j}^2\big(X_{t,j+1}+t\big),\\
&Z_{t,i} = \frac{t}{2q}\big(X_{t,i-1}^2+2X_{t,i}(X_{t,i+1}+t)\big)
\cos\left(\frac1q\sum_{j=1}^q X_{t,j}^2(X_{t,j+1}+t)+X_{t,i}\right).
\end{align*}
We solve this example for $q=2,3,4, 5,$ and  the errors, running time
and convergence rate are shown in Table 4. Again, high order convergence rates are obtained for this coupled example.
\begin{figure}[H]
\centering
\includegraphics[width=3.0in]{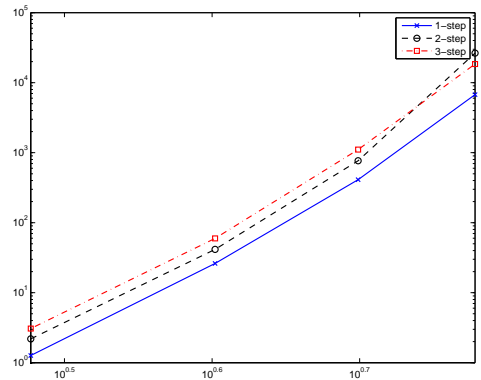}
\caption{Dimensions v.s. the running time
         for Example 2 with $N=128$.}
\end{figure}

\begin{table}
\small
\setlength{\belowcaptionskip}{0pt}
\caption{Errors and convergence rates for Example 3.}
\begin{center}
\begin{tabular}{c|c|c|c|c|c|c|c|c}
\hline
step number    & sparse grid &           &   N=8   &   N=16   &   N=32   &  N=64 &  N=128 & CR         \\
\hline
\multirow{3}{*}{1-step}
& \multirow{3}{*}{$q=2$, $\C_2^3$ \& $\G_2^3$} & $\mathcal{E}_Y$   & 4.246E-03 & 1.925E-03 & 9.181E-04 & 4.508E-04 & 2.243E-04 & 1.058 \\
\cline{3-9}
& \multirow{3}{*}{}         & $\mathcal{E}_Z$ & 1.149E-02 & 4.898E-03 & 2.098E-03 & 9.057E-04 & 3.948E-04 & 1.216 \\
\cline{3-9}
& \multirow{3}{*}{}         & RT & 0.030 & 0.045 & 0.083 & 0.160 & 0.267 &  \\
\cline{3-9}
\hline
\multirow{3}{*}{2-step}
& \multirow{3}{*}{$q=2$, $\C_2^3$ \& $\G_2^3$} & $\mathcal{E}_Y$ & 4.093E-04 & 8.289E-05 & 1.646E-05 & 3.207E-06 & 6.101E-07 & 2.347 \\
\cline{3-9}
& \multirow{3}{*}{}         & $\mathcal{E}_Z$ & 5.429E-03 & 1.395E-03 & 3.568E-04 & 9.099E-05 & 2.315E-05 & 1.968 \\
\cline{3-9}
& \multirow{3}{*}{}         & RT & 0.032 & 0.059 & 0.090 & 0.171 & 0.327 &  \\
\cline{3-9}
\hline
\multirow{3}{*}{3-step}
& \multirow{3}{*}{$q=2$, $\C_2^3$ \& $\G_2^3$} & $\mathcal{E}_Y$ & 2.459E-05 & 2.651E-06 & 2.695E-07 & 2.658E-08 & 2.538E-09 & 3.312 \\
\cline{3-9}
& \multirow{3}{*}{}         & $\mathcal{E}_Z$ & 5.266E-05 & 5.713E-06 & 6.099E-07 & 6.276E-08 & 6.215E-09 & 3.261 \\
\cline{3-9}
& \multirow{3}{*}{}         & RT & 0.033 & 0.075 & 0.132 & 0.276 & 0.430 &  \\
\cline{3-9}
\hline
\multirow{3}{*}{1-step}
& \multirow{3}{*}{$q=3$, $\C_3^4$ \& $\G_3^4$} & $\mathcal{E}_Y$   & 2.463E-03 & 1.068E-03 & 4.812E-04 & 2.229E-04 & 1.072E-04 & 1.130 \\
\cline{3-9}
& \multirow{3}{*}{}         & $\mathcal{E}_Z$ & 3.243E-03 & 1.363E-03 & 5.990E-04 & 2.630E-04 & 1.154E-04 & 1.200 \\
\cline{3-9}
& \multirow{3}{*}{}         & RT & 0.295 & 0.422 & 0.762 & 1.420 & 2.690 &  \\
\cline{3-9}
\hline
\multirow{3}{*}{2-step}
& \multirow{3}{*}{$q=3$, $\C_3^4$ \& $\G_3^4$} & $\mathcal{E}_Y$ & 2.044E-04 & 7.910E-05 & 1.676E-05 & 3.397E-06 & 6.738E-07 & 2.103 \\
\cline{3-9}
& \multirow{3}{*}{}         & $\mathcal{E}_Z$ & 1.073E-03 & 2.444E-04 & 6.285E-05 & 1.632E-05 & 4.227E-06 & 1.988 \\
\cline{3-9}
& \multirow{3}{*}{}         & RT & 0.352 & 0.717 & 1.385 & 2.610 & 4.946 &  \\
\cline{3-9}
\hline
\multirow{3}{*}{3-step}
& \multirow{3}{*}{$q=3$, $\C_3^4$ \& $\G_3^4$} & $\mathcal{E}_Y$ & 8.103E-06 & 1.743E-06 & 1.890E-07 & 1.942E-08 & 1.933E-09 & 3.055 \\
\cline{3-9}
& \multirow{3}{*}{}         & $\mathcal{E}_Z$ & 7.540E-06 & 1.577E-06 & 1.795E-07 & 1.933E-08 & 1.995E-09 & 3.012 \\
\cline{3-9}
& \multirow{3}{*}{}         & RT & 0.364 & 1.223 & 1.936 & 3.743 & 7.150 &  \\
\cline{3-9}
\hline
\multirow{3}{*}{1-step}
& \multirow{3}{*}{$q=4$, $\C_4^5$ \& $\G_4^5$} & $\mathcal{E}_Y$   & 1.663E-03 & 7.351E-04 & 3.068E-04 & 1.370E-04 & 6.294E-05 & 1.187 \\
\cline{3-9}
& \multirow{3}{*}{}         & $\mathcal{E}_Z$ & 1.360E-03 & 5.763E-04 & 2.374E-04 & 1.068E-04 & 4.838E-05 & 1.206 \\
\cline{3-9}
& \multirow{3}{*}{}         & RT & 3.916 & 8.050 & 15.738 & 29.987 & 55.418 &  \\
\cline{3-9}
\hline
\multirow{3}{*}{2-step}
& \multirow{3}{*}{$q=4$, $\C_4^5$ \& $\G_4^5$} & $\mathcal{E}_Y$ & 6.776E-05 & 4.165E-05 & 1.569E-05 & 3.408E-06 & 6.925E-07 & 1.684 \\
\cline{3-9}
& \multirow{3}{*}{}         & $\mathcal{E}_Z$ & 3.665E-04 & 8.033E-05 & 1.727E-05 & 4.479E-06 & 1.191E-06 & 2.070 \\
\cline{3-9}
& \multirow{3}{*}{}         & RT & 5.064 & 13.279 & 28.974 & 56.345 & 105.155 &  \\
\cline{3-9}
\hline
\multirow{3}{*}{3-step}
& \multirow{3}{*}{$q=4$, $\C_4^5$ \& $\G_4^6$} & $\mathcal{E}_Y$ & 4.932E-06 & 1.065E-06 & 1.111E-07 & 1.115E-08 & 1.087E-09 & 3.087 \\
\cline{3-9}
& \multirow{3}{*}{}         & $\mathcal{E}_Z$ & 2.474E-06 & 5.495E-07 & 6.062E-08 & 6.348E-09 & 6.390E-10 & 3.027 \\
\cline{3-9}
& \multirow{3}{*}{}         & RT & 34.015 & 117.047 & 250.853 & 483.481 & 926.454 &  \\
\cline{3-9}
\hline
\multirow{3}{*}{1-step}
& \multirow{3}{*}{$q=5$, $\C_5^6$ \& $\G_5^6$} & $\mathcal{E}_Y$ & 1.167E-03 & 5.489E-04 & 2.341E-04 & 9.343E-05 & 4.184E-05 & 1.216 \\
\cline{3-9}
& \multirow{3}{*}{}         & $\mathcal{E}_Z$ & 6.751E-04 & 3.055E-04 & 1.268E-04 & 5.183E-05 & 2.408E-05 & 1.218 \\
\cline{3-9}
& \multirow{3}{*}{}         & RT & 58.937 & 134.345 & 289.234 & 570.587 & 1049.247 &  \\
\cline{3-9}
\hline
\multirow{3}{*}{2-step}
& \multirow{3}{*}{$q=5$, $\C_5^6$ \& $\G_5^7$} & $\mathcal{E}_Y$ & 6.572E-05 & 4.173E-05 & 1.246E-05 & 2.548E-06 & 5.087E-07 & 1.806 \\
\cline{3-9}
& \multirow{3}{*}{}         & $\mathcal{E}_Z$ & 1.464E-04 & 2.885E-05 & 6.466E-06 & 1.763E-06 & 4.792E-07 & 2.054 \\
\cline{3-9}
& \multirow{3}{*}{}         & RT & 736.405 & 2004.185 & 4124.286 & 7776.645 & 14561.700 &  \\
\cline{3-9}
\hline
\multirow{3}{*}{3-step}
& \multirow{3}{*}{$q=5$, $\C_5^6$ \& $\G_5^7$} & $\mathcal{E}_Y$ & 1.649E-06 & 6.334E-07 & 9.009E-08 & 9.209E-09 & 9.147E-10 & 2.774 \\
\cline{3-9}
& \multirow{3}{*}{}         & $\mathcal{E}_Z$ & 5.318E-07 & 2.072E-07 & 3.079E-08 & 3.274E-09 & 3.372E-10 & 2.723 \\
\cline{3-9}
& \multirow{3}{*}{}         & RT & 709.587 & 2599.783 & 5916.660 & 11366.633 & 21380.402 &  \\
\cline{3-9}
\hline
\end{tabular}
\end{center}
\end{table}

\section{Conclusions}

In this work, we have extended our previous work \cite{ZFZ2014} of multi-step schemes to solve high-dimensional FBSDEs, by combining the sparse grid spatial discretizations and the sparse grid quadrature \& interpolations, and the entire algorithm admits efficient and high accurate approximations in high-dimensions. It is shown that the proposed numerical techniques in this work are more efficient than our previous work in \cite{ZFZ2014}. However, we remark that how to build efficient high dimensional approximations is a long-term open question, and more efforts are still needed. Possible extensions along this direction include the adaptive sparse grid approaches, the anisotropic constructions of spatial approximations. Also, a rigorous error analysis will also be part of our future studies.

\end{document}